\documentstyle[11pt,amstex]{amsart}
\newcount\skewfactor
\def\mathunderaccent#1#2 {\let\theaccent#1\skewfactor#2
\mathpalette\putaccentunder}
\def\putaccentunder#1#2{\oalign{$#1#2$\crcr\hidewidth
\vbox to.2ex{\hbox{$#1\skew\skewfactor\theaccent{}$}\vss}\hidewidth}}
\def\name{\mathunderaccent\tilde-3 }

\newcommand{\cf}{{\rm cf}}
\newcommand{\otp}{{\rm otp}}
\newcommand{\dom}{\operatorname {dom}} 
\newcommand{\rng}{\operatorname {rng}}
\newcommand{\con}{{\frak c}}
\newcommand{\rest}{{\,|\grave{}\,}}
\newcommand{\lesdot}{\mathrel{\mathord{<}\!\!\raise 0.8
pt\hbox{$\scriptstyle\circ$}}}  
\newcommand{\dk}{{{\rm d}_\kappa}}
\newcommand{\gek}{{{\frak e}_\kappa}}
\newcommand{\gdk}{{{\frak d}_\kappa}}
\newcommand{\gbk}{{{\frak b}_\kappa}}
\newcommand{\gec}{{{\frak e}_\con}}
\newcommand{\dc}{{{\rm d}_\con}}
\newcommand{\Stk}{{\mathcal S}^\theta_\kappa}
\newcommand{\ctk}{{\rm c}(\kappa,\theta)}
\newcommand{\ctkm}{{\rm c}^-(\kappa,\theta)}
\newcommand{\btk}{{\rm b}(\kappa,\theta)}
\newcommand{\btkm}{{\rm b}^-(\kappa,\theta)}
\newcommand{\dks}{{{\rm d}^*_\kappa}}
\newcommand{\geks}{{{\frak e}^*_\kappa}}
\newcommand{\gecs}{{{\frak e}^*_\con}}
\newcommand{\dcs}{{{\rm d}^*_\con}}
\newcommand{\kk}{{}^{\textstyle \kappa}\kappa}
\newcommand{\bP}{{\mathbb P}}
\newcommand{\bQ}{{\mathbb Q}}
\newcommand{\bR}{{\mathbb R}}
\newcommand{\bV}{{\bf V}}

\newcommand{\lh}{\ell g\/} 
\newcommand{\tcf}{{\rm tcf\/}}
\newcommand{\jbt}{{J^{\rm bd}_\theta}}
\newcommand{\jbd}{{J^{\rm bd}_\delta}}
\newcommand{\jbk}{{J^{\rm bd}_\kappa}}
\newcommand{\cD}{{\mathcal D}}
\newcommand{\cA}{{\mathcal A}}
\newcommand{\cI}{{\mathcal I}}
\newcommand{\diagun}{\mathop{\triangle}}
\newcommand{\APf}{{\mathcal AP}^f}
\newcommand{\forces}{\mathrel {{\vrule height 6.9pt depth -0.1pt}\!\vdash}}

\newcommand{\QED}{\hfill\vrule width 6pt height 6pt depth 0pt\vspace{0.1in}} 
\newcommand{\Proof}{\noindent{\sc Proof} \hspace{0.2in}} 

\newtheorem{theorem}{Theorem}[section] 
\newtheorem{claim}{Claim}[theorem]
\newtheorem{proposition}[theorem]{Proposition} 
\newtheorem{lemma}[theorem]{Lemma} 

\theoremstyle{definition}
\newtheorem{definition}[theorem]{Definition}
\newtheorem{notation}[theorem]{Notation}

\theoremstyle{remark}
\newtheorem{remark}[theorem]{Remark}
\newtheorem{conclusion}[theorem]{Conclusion}

\title[\it On Ciesielski's Problems]{
\vspace {3.0cm}\uppercase 
{\Large \bf On Ciesielski's Problems}} 
\author[S. Shelah]{\uppercase {\bf S. Shelah}}
\address{Institute of Mathematics\\
The Hebrew University\\
Jerusalem 91904, Israel\\
and Rutgers University\\
Mathematics Department\\
New Brunswick, NJ 08854, USA
}
\email{shelah@@math.huji.ac.il}
\date{\today} 

\thanks{The research was partially supported by ``Basic Research Foundation''
of the Israel Academy of Sciences and Humanities. Publication 675.}

\subjclass{Primary: 26Axx, 03E05;\\ Secondary 03E35, 03E50}

\begin{document}

\maketitle 

\bigskip
\bigskip
\bigskip
{\it Abstract.}
In the present paper we discuss some problems formulated in Ciesielski
\cite{Ci98}.
\bigskip
\bigskip
\bigskip
\setcounter{section}{-1}

\stepcounter{section}
\subsection*{\quad 0. Introduction}
I was asked to read and comment on Ciesielski's survey paper \cite{Ci98}. I
have found it very exciting and illuminating. Quite naturally I was not able
to resist the temptation to look mainly at the open problems formulated in
this nice paper. Some of them are related to my research {\em in progress} and
may be solved soon. This is in particular the case with Problems 5, 1 and 6,
for which relevant information should be given by Ciesielski Shelah
\cite{Sh:F288} and Ros{\l}anowski Shelah \cite{RoSh:670}. For some other
problems (like \cite[Problem 3]{Ci98}) I have ideas that could work and this
may materialize in a continuation of the present paper. 

Here we would like to present answers to three problems and address a fourth
one. In the first section we solve \cite[Problem 8]{Ci98} and we show that,
consistently, $\dc$ is a singular cardinal and $\gec<\dc$ (in \ref{conjeden};
see \ref{cardinv} for the definitions of $\gec,\dc$). In the next section we
present some  results relevant for \cite[Problem 9]{Ci98}. We do not solve the
problem, but it was formulated in a very general way ({\em When does
$\dc=\dcs$ or $\gec=\gecs$ hold?}) making the full answer rather difficult. 
The third section answers \cite[Problem 7]{Ci98}. We show there that the
Martin Axiom for $\sigma$--centered forcing notions implies that for every
function $f:\bR^2\longrightarrow\bR$ there are functions $g_n,h_n:\bR
\longrightarrow\bR$, $n<\omega$, such that $f(x,y)=\sum\limits_{n=0}^{\infty}
g_n(x)\cdot h_n(y)$. Finally, in the next section we deal with countably
continuous functions and we show (in \ref{pr4thm}) that in the Cohen model
they are exactly the functions $f$ with the property that
\[(\forall U\in [\bR]^{\textstyle \aleph_1})(\exists U^*\in [U]^{\textstyle
\aleph_1})(f\rest U^*\mbox{ is continuous\/}).\]
This answers negatively \cite[Problem 4]{Ci98}.

\noindent{\bf Notation:} Our notation is rather standard and compatible with
that of classical textbooks on Set Theory (like Jech \cite{J} or
Bartoszy\'nski Judah \cite{BaJu95}). However in forcing we keep the convention
that {\em a stronger condition is the larger one}. 

\medskip
\begin{notation}
\label{notacja}
We will keep the following rules for our notation:
\begin{enumerate}
\item $\alpha,\beta,\gamma,\delta,\xi,\zeta, i,j\ldots$ will denote ordinals.
\item $\kappa,\lambda,\mu,\theta\ldots$ will stand for cardinal numbers,
$\con$ is the cardinality of the continuum.
\item A bar above a name indicates that the object is a sequence, usually
$\bar{X}$ will be $\langle X_i: i<\lh(\bar{X})\rangle$, where $\lh(\bar{X})$
denotes the length of $\bar{X}$.
\item A tilde indicates that we are dealing with a name for an object in
forcing extension (like $\name{x}$). The canonical $\bP$--name for a generic
filter is called $\name{G}_{\bP}$.
\item For a cardinal $\kappa$, the quantifiers $(\exists^\kappa i)$ and
$(\forall^\kappa i)$ are abbreviations for ``there is $\kappa$ many $i$ such
that$\ldots$'' and ``for all but less than $\kappa$ many $i\ldots$'', 
respectively.
\item $\otp$ stands for ``order type''. When using elements of the pcf--theory
we will follow the notation and terminology of \cite{Sh:g}. In particular,
$\tcf$ will stand for ``true cofinality'' and $\jbt$ will denote the ideal of
bounded subsets of $\theta$.
\end{enumerate}
\end{notation}

\stepcounter{section}
\subsection*{\quad 1. Around $\dk$ and $\gek$}
\begin{definition}
\label{cardinv}
Let $\theta\leq\kappa$ be cardinals. 
\begin{enumerate}
\item Let $\Stk\stackrel{\rm def}{=}\prod\limits_{i<\kappa}[\kappa]^{
\textstyle <\!\theta}$. 
\item We define the following cardinal coefficients of the space $\kk$:
\[\begin{array}{c}
\dk =\min\{|F|:F\subseteq\kk\ \&\ (\forall g\in\kk)(\exists f\in F)(
\exists^\kappa i\!<\! \kappa)(f(i)=g(i))\},\\
\gek=\min\{|F|:F\subseteq\kk\ \&\ (\forall g\in\kk)(\exists f\in F)(
\forall^\kappa i\!<\! \kappa)(f(i)\neq g(i))\},\\
\gdk=\min\{|F|:F\subseteq\kk\ \&\ (\forall g\in\kk)(\exists f\in F)(
\forall^\kappa i\!<\! \kappa)(g(i)<f(i))\},\\
\gbk=\min\{|F|:F\subseteq\kk\ \&\ (\forall g\in\kk)(\exists f\in F)(
\exists^\kappa i\!<\! \kappa)(g(i)<f(i))\},\\
\ctk=\min\{|G|:G\subseteq\Stk\ \&\ (\forall g\in\kk)(\exists\bar{S}\in G)(
\forall^\kappa i\!<\! \kappa)(g(i)\in S_i)\},\\
\ctkm=\min\{|G|:G\subseteq\Stk\ \&\ (\forall g\in\kk)(\exists\bar{S}\in G)(
\exists^\kappa i\!<\! \kappa)(g(i)\in S_i)\},\\
\btk=\min\{|F|:F\subseteq\kk\ \&\ (\forall\bar{S}\in\Stk)(\exists f\in F)(
\forall^\kappa i\!<\! \kappa)(f(i)\notin S_i)\}.\\
\btkm=\min\{|F|:F\subseteq\kk\ \&\ (\forall\bar{S}\in\Stk)(\exists f\in F)(
\exists^\kappa i\!<\! \kappa)(f(i)\notin S_i)\}.\\
  \end{array}\]
\item For functions $f,g\in\kk$ we say that {\em $f$ dominates $g$} (in short:
$g<^*_\kappa f$) if $(\forall^\kappa i\!<\!\kappa)(g(i)<f(i))$.\\
\relax [Thus $\gbk$ and $\gdk$ are the unbounded number and the dominating
number, respectively, of the partial order $(\kk,{<^*_\kappa})$.]
\end{enumerate}
\end{definition}

\begin{remark}
\label{rempro}
\begin{enumerate}
\item The cardinal invariants introduced in \ref{cardinv} are natural
generalizations of those studied in Set Theory of the Reals; see e.g.\
Bartoszy\'nski Judah \cite{BaJu95} or Goldstern Shelah \cite{GoSh:448}.  
\item Using \ref{cardinv}, we may reformulate \cite[Problem 8]{Ci98} as
follows: 
\begin{enumerate}
\item[(a)] Is it consistent that $\dc>\gec$?
\item[(b)] Can $\dc$ be a singular cardinal?
\end{enumerate}
(see \cite[4.7, 4.12]{Ci98}).
\end{enumerate}
\end{remark}

\begin{proposition}
\label{trivinv}
\begin{enumerate}
\item The partial order $(\kk,{<^*_\kappa})$ is $\gbk$--directed. The cardinal
$\gbk$ is regular. If $\kappa$ is regular then $\gbk={\rm
b}^-(\kappa,\kappa)$.  
\item $\gbk\leq\dk$. If $\kappa$ is a successor then $\dk=\gbk$.
\item $\cf(\kappa)<{\rm c}^-(\kappa,\kappa)$ and $\theta<\kappa\ \Rightarrow\
\kappa<c^-(\kappa,\theta^+)$. 
\item Assume that either  $\theta<\cf(\kappa)$ or $\theta=\cf(\kappa)$ is a
successor cardinal. Then $\dk=\ctkm$.
\end{enumerate}
\end{proposition}

\Proof 1) and 3)\quad Should be clear.
\medskip

\noindent 2)\quad For a function $f\in\kk$ let $f^+\in\kk$ be defined by
$f^+(i)=f(i)+1$. Clearly, if $F\subseteq\kk$ is a family witnessing the
minimum in the definition of $\dk$ then $\{f^+:f\in F\}$ is a
$<^*_\kappa$--unbounded family. Hence $\gbk\leq\dk$.

Assume now that $\kappa=\mu^+$ and let $F\subseteq\kk$ be
$<^*_\kappa$--unbounded, $|F|=\gbk$. Note that necessarily $\gbk>\kappa$. For
each $\alpha<\kappa$ fix a sequence $\langle\beta_{\alpha,\xi}:\xi<\mu\rangle$
such that $\{\beta_{\alpha,\xi}:\xi<\mu\}=\alpha+1$. For $f\in F$ and $\xi<
\mu$ let $h^f_\xi\in\kk$ be such that $(\forall i<\kappa)(h^f_\xi(i)=\beta_{f
(i),\xi})$. Let 
\[F^*\stackrel{\rm def}{=}\{h^f_\xi: f\in F\ \&\ \xi<\mu\}.\]
Then $|F^*|\leq |F|+\mu=\gbk$. Suppose $g\in\kk$. By the choice of $F$, we
find $f\in F$ such that the set $A\stackrel{\rm def}{=}\{i<\kappa:
g(i)<f(i)\}$ is of cardinality $\kappa$. For $i\in A$ let $\xi_i<\mu$ be such
that $g(i)=\beta_{f(i),\xi_i}$. Then for some $\xi<\mu$ the set $A_\xi=\{i\in
A: \xi_i=\xi\}$ is of size $\kappa$. Look at the function $h^f_\xi\in F^*$:
for every $i\in A_\xi$ we have $g(i)=h^f_\xi(i)$. 
\medskip

\noindent 4)\quad First note that plainly $\ctkm\leq\dk$, so we have to show
the converse inequality (under our assumptions).

Assume $\theta<\cf(\kappa)$. Let $G\subseteq\Stk$ be such that $|G|=\ctkm$ and 
\[(\forall g\in\kk)(\exists\bar{S}\in G)(\exists^\kappa i\!<\! \kappa)(g(i)
\in S_i).\]
For $\bar{S}\in G$ and $i<\kappa$ fix an enumeration $\{\beta^{\bar{S},
i}_\varepsilon:\varepsilon<\varepsilon^{\bar{S},i}\}$ of $S_i$ (so
$\varepsilon^{\bar{S},i}<\theta$).  
Next define functions $h^{\bar{S}}_\varepsilon\in\kappa$ (for $\bar{S}\in G$
and $\varepsilon<\theta$) by
\[h^{\bar{S}}_\varepsilon(i)=\left\{
\begin{array}{ll}
\beta^{\bar{S},i}_\varepsilon&\mbox{ if }\varepsilon<\varepsilon^{\bar{S},
i}_\varepsilon\\
0  &\mbox{ otherwise.}
\end{array}\right.\]
Suppose that $g\in\kk$. Take $\bar{S}\in G$ such that $(\exists^\kappa i\!<\!
\kappa)(g(i)\in S_i)$. Then for some $\varepsilon<\theta$ we have 
\[(\exists^\kappa i\!<\!\kappa)(g(i)=\beta^{\bar{S},i}_\varepsilon=
h^{\bar{S}}_\varepsilon(i)),\]
and hence we may conclude that $\dk\leq \ctkm+\theta$ is witnessed by the
family $\{h^{\bar{S}}_\varepsilon: \bar{S}\in G\ \&\ \varepsilon<\theta\}$. 
Finally we note that $\ctkm+\theta=\ctkm$ (by (3); remember $\ctkm\geq {\rm
c}^-(\kappa,\kappa)$).

If $\theta=\cf(\kappa)$ is a successor cardinal, say $\theta=\mu^+$, then we
proceed similarly: we may assume that for each $\bar{S}\in G$ and $i<\kappa$
we have $|S_i|=\mu=\varepsilon^{\bar{S},i}$ and we finish as above (as
$\mu<\cf(\kappa)$). \QED

\begin{proposition}
\label{getdk}
Assume that $\kappa$ is a strong limit singular cardinal, $\cf(\kappa)=\theta
>\aleph_0$. Then 
\[\ctkm=\ctk=\dk=2^\kappa.\]
\end{proposition}

\Proof Clearly $\kappa<\ctkm\leq\ctk\leq 2^\kappa$ (remember
\ref{trivinv}(3)) and $\ctkm\leq\dk\leq 2^\kappa$, so it suffices to show that
$\ctkm\geq 2^\kappa$.  

Suppose that $G\subseteq \Stk$, $|G|=\mu$, $\kappa<\mu<\mu^+\leq 2^\kappa$.

Choose an increasing continuous sequence $\langle\kappa_i:i<\theta\rangle$
such that 
\[\theta<\kappa_0\quad\mbox{ and }\quad\sup\limits_{i<\theta}\kappa_i=\kappa
\quad\mbox{ and }\quad\big(\forall i<\theta)(2^{\sum\limits_{j<i}\kappa_j+
\aleph_0}<\kappa_i\big).\]
Next, using \cite[Ch. VIII, \S1]{Sh:g}, pick $\bar{\chi}=\langle\chi_i: i<
\theta\rangle$ such that
\begin{enumerate}
\item[(i)]   $\bar{\chi}$ is a strictly increasing sequence of regular
cardinals,
\item[(ii)]  $\kappa_i<\chi_i<\kappa$ for each $i<\theta$,
\item[(iii)] $\tcf(\prod\limits_{i<\theta}\chi_i/\jbt)=\mu^+$.
\end{enumerate}
Now, for every $\bar{S}\in G$ define a function $h^{\bar{S}}\in\prod
\limits_{i<\theta}\chi_i$ by
\[h^{\bar{S}}(i)=\sup\{\alpha<\chi_i: \alpha\in S_\gamma\ \mbox{ and }\
\gamma<\kappa_i\}.\]
Note that $|\{\alpha<\chi_i: \alpha\in S_\gamma,\ \gamma<\kappa_i\}|\leq
\kappa_i\cdot\theta<\chi_i$, so (as $\chi_i$ is regular) $h^{\bar{S}}(i)<
\chi_i$. It follows from (iii) that there a function $h\in\prod\limits_{i<
\theta}\chi_i$ such that 
\[(\forall \bar{S}\in G)(h^{\bar{S}}<_{\jbt} h).\] 
Finally define a function $g\in\kk$ by:
\[\mbox{if }\quad\sup\limits_{j<i}\kappa_j\leq\gamma<\kappa_i\quad\mbox{ then
}\quad g(\gamma)=h(i).\]
Note that for each $\bar{S}\in G$ we have
\[\{\gamma<\kappa: g(\gamma)\in S_\gamma\}\subseteq \bigcup\big\{[\sup_{j<i}
\kappa_j,\kappa_i): i<\theta,\ h^{\bar{S}}(i)\geq h(i)\big\}\subseteq
\kappa_{j(\bar{S})},\]
where $j(\bar{S})=\min\{j<\theta:\{i<\theta:h^{\bar{S}}(i)\geq h(i)\}\subseteq
j\}$. Consequently the function $g$ shows that the family $G$ cannot
witness the minimum in the definition of $\ctkm$ and we are done. \QED

\begin{remark}
Actually much weaker assumptions are sufficient to get the conclusion of
\ref{getdk}. For example, almost always we may allow $\theta=\aleph_0$ (see
\cite{E12}). 
\end{remark}

\begin{proposition}
\label{getek}
If $\kappa$ is a singular cardinal, $\theta<\kappa$ then $\gek=\btk=\kappa^+$.
\end{proposition}

\Proof First note that $\kappa<\gek\leq\btk$, so it is enough to show that
$\btk\leq\kappa^+$. 

By \cite[Ch. II, 1.5]{Sh:g}, we may find an increasing sequence
$\langle\chi_i: i<\cf(\kappa)\rangle$ of regular cardinals cofinal in $\kappa$
and such that 
\[\theta<\chi_0,\quad \tcf(\prod_{i<\cf(\kappa)}\chi_i/J^{{\rm bd}}_{\cf(
\kappa)})=\kappa^+\quad\mbox{and}\quad(\forall i<\cf(\kappa))(\sup_{j<i}\chi_j
<\chi_i<\kappa).\] 
Let $\langle h_\alpha:\alpha<\kappa^+\rangle\subseteq\prod\limits_{i<\cf(
\kappa)}\chi_i$ be a $<_{J^{\rm bd}_{\cf(\kappa)}}$--increasing sequence
cofinal in $(\prod\limits_{i<\cf(\kappa)}\chi_i,<_{J^{\rm bd}_{\cf(\kappa
)}})$. For $i<\cf(\kappa)$ put $\mu_i\stackrel{\rm def}{=}\sup\limits_{j<i}
\chi_j$. Then the sequence $\langle\mu_i:i<\cf(\kappa)\rangle$ is increasing
continuous with limit $\kappa$. Now we define functions $f_\alpha\in\kk$ (for
$\alpha<\kappa^+$) by:
\[\mu_i\leq\xi<\mu_{i+1}\ \&\ i<\cf(\kappa)\quad\Rightarrow\quad f_\alpha(\xi)
=h_\alpha(i+1).\]
We claim that
\[(\forall\bar{S}\in\Stk)(\exists \alpha<\kappa^+)(\forall^\kappa \xi\!<\!
\kappa)(f_\alpha(\xi)\notin S_\xi).\]
So suppose $\bar{S}\in\Stk$. Define a function $h^{\bar{S}}\in\prod
\limits_{i<\cf(\kappa)}\chi_i$ by
\[h^{\bar{S}}(i)=\sup\{\alpha<\chi_i: \alpha\in S_\xi\ \mbox{ and }\
\xi<\mu_i\}\]
(note that the set on the right-hand side of the formula above is of size
$<\chi_i$ so the supremum is below $\chi_i$). Take $\alpha<\kappa^+$ and
$j^*<\cf(\kappa)$ such that 
\[j^*\leq i<\cf(\kappa)\quad\Rightarrow\quad h^{\bar{S}}(i)<h_\alpha(i),\]
and note that then
\[\{\xi<\kappa: f_\alpha(\xi)\in S_\xi\}\subseteq \mu_{j^*}.\]
So we are done. \QED

\begin{proposition}
\label{sah1}
If $\kappa$ is singular and $\theta<\kappa$ then
\begin{enumerate}
\item[(a)] $\dk\geq\ctkm\geq {\rm pp}_{J^{\rm bd}_{\cf(\kappa)}}(\kappa)$,
\item[(b)] $\dk\geq\gek$ and if $\dk=\gek$ then
\[\cf(\kappa)\leq\theta<\kappa\quad\Rightarrow\quad {\rm pp}_\theta(\kappa)
=\kappa^+.\]
\end{enumerate}
\end{proposition}

\Proof To show clause (a) take any $\mu<{\rm pp}_{J^{\rm bd}_{\cf(\kappa)}}(
\kappa)$ and essentially repeat the proof of \ref{getek} for $\mu^+$ (remember
\cite[Ch. II, 2.3]{Sh:g}). The assertion (b) follows from (a) and
\ref{getek}. \QED 

\begin{proposition}
\label{forinv}
Assume that $\bP$ is a $\cf(\theta)$--cc forcing notion.
\begin{enumerate}
\item $\forces_{\bP}$`` $(\forall \bar{S}\in\Stk)(\exists\bar{S}^*\in\Stk\cap
\bV)(\forall i<\kappa)(S_i\subseteq S_i^*)$ ''. 
\item $\forces_{\bP}$`` $\ctk\geq (\ctk)^\bV$ and $\ctkm\geq (\ctkm)^\bV$ ''. 
\item If $\kappa^{<\cf(\theta)}=\kappa$ then 
\[\forces_{\bP}\mbox{`` }\ctk=(\ctk)^\bV\mbox{ and }\ctkm=(\ctkm)^\bV\mbox{
''.}\]   
\item If $\theta=\cf(\theta)<\kappa$ and either $\theta<\cf(\kappa)$ or
$\theta=\cf(\kappa)$ is a successor cardinal then $\forces_{\bP}$``
$\dk=(\dk)^\bV$ ''.  
\end{enumerate}
\end{proposition}

\Proof 1)\quad Suppose that $\name{A}$ is a $\bP$--name for a set of ordinals,
$\forces_{\bP}|\name{A}|<\theta$. Since $\bP$ satisfies the $\cf(\theta)$--cc,
we find a cardinal $\mu<\theta$ and a $\bP$--name $\name{h}$ such that
$\forces_{\bP}$`` $\name{h}:\mu\stackrel{\rm onto}{\longrightarrow}
\name{A}$ ''. By the $\cf(\theta)$--cc again, we find sets $B_i$ (for $i<\mu$)
such that $|B_i|<\cf(\theta)$ and $\forces_{\bP}\name{h}(i)\in B_i$. Let
$A=\bigcup\limits_{i<\mu} B_i$. Then $\forces_{\bP}\name{A}\subseteq A$ and:
if $\cf(\theta)<\theta$ then $|A|\leq\mu\cdot\cf(\theta)<\theta$ and if
$\cf(\theta)=\theta$ then $|A|<\theta$ as $\mu<\cf(\theta)$. The rest should
be clear. 
\medskip

\noindent 4)\quad Let $F\subseteq\kk$, $F\in\bV$ be a family witnessing
the minimum in the definition of $\dk$. We are going to show that
\[\forces_{\bP}\mbox{`` }(\forall g\in\kk)(\exists f\in F)(\exists^\kappa i\!
<\!\kappa)(g(i)=f(i))\mbox{ ''.}\]
So suppose that $p\in\bP$ and $\name{g}$ are such that $p\forces$``$\name{g}
\in\kk$''. Choose a sequence $\langle p_i: i<\kappa\rangle$ of conditions and
a function $g\in\kk$ such that 
\[(\forall i<\kappa)(p\leq p_i\ \ \&\ \ p_i\forces_{\bP}\name{g}(i)=g(i)).\]
By the choice of $F$ we find $f\in F$ such that the set $A\stackrel{\rm
def}{=}\{i<\kappa: g(i)=f(i)\}$ is of size $\kappa$. Next choose a condition
$q\geq p$ such that
\[q\forces_{\bP}\mbox{`` }|\{i\in A: p_i\in\name{G}_{\bP}\}|=\kappa\mbox{
''}.\]
[Possible, as otherwise $p\forces$``$|\{i\in A: p_i\in\name{G}_{\bP}\}|\leq
\mu$'' for some $\mu<\kappa$ (remember that $\theta\leq\cf(\kappa)$). So we
have a $\bP$--name $\name{h}$ for a function from $\mu$ into $A$ such that
$p\forces (\forall i\in A)(p_i\in\name{G}_{\bP}\ \Rightarrow\ i\in\rng(
\name{h}))$. For each $\zeta\in\mu$ the set $B_\zeta=\{i\in A: (\exists p'\geq
p)(p'\forces \name{h}(\zeta)=i)\}$ is of size $<\theta$ and hence
$|\bigcup\limits_{\zeta<\mu} B_\zeta|\leq\theta\cdot\mu<\kappa$. Take any
$i\in A\setminus\bigcup\limits_{\zeta<\mu} B_\zeta$ and look at the condition
$p_i$.]\\
Now note that the condition $q$ forces `` $(\exists^\kappa i\in A)(\name{g}(i)
=g(i))$ ''. 

Thus we have proved that $\forces_{\bP}$``$\dk\leq(\dk)^{\bV}$''. For the
converse inequality we use \ref{trivinv}(4) and \ref{forinv}(2). Thus we get
\[\forces_{\bP}\mbox{`` }\dk=\ctkm\geq (\ctkm)^{\bV}=(\dk)^{\bV}\mbox{ '',}\]
finishing the proof. \QED
\medskip

Now may get the affirmative answer to \cite[Problem 8]{Ci98} (see
\ref{rempro}(2)): 

\begin{conclusion}
\label{conjeden}
It is consistent that $\dc$ is a singular cardinal and $\gec<\dc$ (modulo
existence of high enough measurables).
\end{conclusion}

\Proof First we force that there is $\kappa$ satisfying the assumptions of
\ref{getdk} and such that $2^\kappa$ singular. How? Start with a supercompact
Laver indestructible $\kappa$ and make $2^\kappa$ to have cofinality
$\kappa^+$, $\kappa$ still supercompact. Next force $\kappa$ to have
cofinality $\aleph_1$, say as in Magidor \cite{Mg4}. (By \cite{Sh:137} we can
make $\kappa$ to be the $\omega_1$--th fix point among the alephs.) So now we
have $\dk=2^\kappa$, $\cf(2^\kappa)=\kappa^+<2^\kappa$. Next add $\kappa$
Cohen reals. Since this forcing satisfies the $\aleph_1$--cc and is of
cardinality $\kappa$ we conclude that, by \ref{forinv}(4), in the final
universe $\dk=\dc$ remains the same (so it is singular). Finally, by
\ref{getek}, we know that in the resulting model $\gec=\kappa^+<\dc$. \QED

\begin{remark}
\begin{enumerate}
\item In fact, if we waive the requirement ``$\dc$ is singular'' then
$2^\kappa=\kappa^{++}$ is enough for the proof, so we can get even
$\kappa=\aleph_{\omega_1}$, $\dk=\aleph_{\omega_1+2}$ and $\gek=\aleph_{
\omega_1+1}$.
\item What is the consistency strength? By Gitik \cite{Gi1} the consistency
strength of
\begin{enumerate}
\item[$(\oplus)$] \qquad $\kappa$ is measurable and $2^\kappa>\kappa^+$ 
\end{enumerate}
is that of the existence of a measurable cardinal $\kappa$ of Mitchell order
$\kappa^{++}$. By Gitik \cite{Gi2}, the consistency strength of
\begin{enumerate}
\item[$(\otimes)$] \qquad $\kappa$ is measurable and $2^\kappa=\aleph_{
\kappa^+}$ 
\end{enumerate}
is that of the existence of a hypermeasurable cardinal with sequence of
measures of length $\aleph_{\kappa^+}$. 
\end{enumerate}
\end{remark}

\stepcounter{section}
\subsection*{\quad 2. Around $\dks$ and $\geks$}
In this section we address \cite[Problem 9]{Ci98}. The problem reads
\begin{quotation}
{\em 
When does $\dc=\dcs$ or $\gec=\gecs$ hold?
}
\end{quotation}
(see \ref{invstar} for the definitions of $\dcs$, $\gecs$). Though we do not
answer the question fully, we are able to give examples of situations in which
the equalities hold. The results here should be combined with those from the
previous section, of course.

\begin{definition}
\label{invstar}
We define the following cardinal coefficients of the space $\kk$:
\[\begin{array}{l}
\dks=\\
\min\{|F|:F\subseteq\kk\ \&\ (\forall G\!\in\! [\kk]^{\textstyle\kappa})(
\exists f\!\in\! F)(\forall g\!\in\! G)(\exists^\kappa i\!<\!\kappa)(f(i)
=g(i))\},\\  
\ \\
\geks=\\
\min\{|F|:F\subseteq\kk\ \&\ (\forall G\!\in\! [\kk]^{\textstyle\kappa})(
\exists f\!\in\! F)(\forall g\!\in\! G)(\forall^\kappa i\!<\!\kappa)(f(i)\neq
g(i))\},\\ 
  \end{array}\]
\end{definition}

\begin{proposition}
\label{notmotiv}
\begin{enumerate}
\item \mbox{\rm [Ciesielski and Jordan; see Jordan \cite{Jo}]} If
$\kappa=\kappa^{<\kappa}$ then $\dk=\dks$ and $\gek=\geks$. 
\item If $\kappa$ is a successor cardinal then $\dks\leq \cf([\dk]^{\textstyle
\kappa},\subseteq)$.  
\item Suppose that $\bar{\lambda}$ is an increasing sequence cofinal in
$\kappa$, $\lh(\bar{\lambda})=\delta\leq\kappa$ such that $\tcf\big(
\prod\limits_{i<\delta}\lambda_i/\jbd\big)=\theta$. Then $\gek\leq\geks\leq
\theta\leq\dk\leq\dks$.
\item If $\kappa$ is singular then $\gek=\geks=\kappa^+$.
\end{enumerate}
\end{proposition}

\Proof 3)\quad Repeat the proof of \ref{getek}, \ref{sah1} with suitable
(minor) changes, see below too. 
\medskip

\noindent  4)\quad A minor modification of the proof of \ref{getek} shows
it. Proceed like there, but instead of functions $h^{\bar{S}}$ (for $\bar{S}
\in\Stk$) consider functions $h_{\bar{g}}\in\prod\limits_{i<\cf(\kappa)}
\chi_i$ (for $\bar{g}=\langle g_\xi:\xi<\kappa\rangle\subseteq\kk$) defined by
\[h_{\bar{g}}(i)=\sup\{\alpha<\chi_i:\alpha=g_\xi(\zeta),\
\xi,\zeta<\mu_i\}.\]
The rest should be clear. \QED

\begin{remark}
\label{remwhat}
\begin{enumerate}
\item Concerning the assumptions of \ref{fildkek} below, note that by Gitik
Shelah \cite{GiSh:597} there may be such ultrafilters in various cases. 
Necessarily $\theta\geq\kappa^+$; it can be $\kappa^+$, which is the
interesting case, and can have $2^\kappa$ singular. See more in D\v{z}amonja
Shelah \cite{DjSh:659}.
\item Concerning \ref{fildkek}(c), note that {\em if\/} $\kappa$ is just
strongly inaccessible, $\kappa<\mu$ and $f\in\kk$ satisfies $(\forall
i<\kappa)(\cf(f(i))>|i|)$ (more if we want to preserve being a large cardinal)
{\em then\/} there is a $\kappa$--strategically closed $\kappa^+$--cc forcing
notion $\bQ$ such that $\forces_{\bQ}$`` $\prod\limits_{i<\kappa} f(i)/\jbk$
is $\mu$--directed ''. (Just iterate the forcing adding $g\in\prod\limits_{
i<\kappa}f(i)$ dominating all members of $\prod\limits_{i<\kappa}f(i)$ from
the ground model; so a condition fixes $g\rest \alpha$ (for some
$\alpha<\kappa$) and promises $g \geq^* g_0\in\prod\limits_{i<\kappa}
f(i)$.)
\item So under the assumption of \ref{fildkek} (all parts), $\geks\leq\dk$.
\end{enumerate}
\end{remark}

\begin{proposition}
\label{fildkek}
Suppose that cardinals $\theta,\kappa$ are such that there is a normal
ultrafilter $\cD$ on $\kappa$ generated by $\theta$ sets. Then
\begin{enumerate}
\item[(a)] $\geks\leq\theta+\kappa^+$,
\item[(b)] {\em if}\ for every family $\cA\subseteq\cD$ of size $<\mu$ there
is $B\in \cD$ such that $(\forall A\in\cA)(B\subseteq^* A)$ {\em then}\
$\mu\leq\dk$, 
\item[(c)] {\em if}\ there is a function $f\in\kk$ such that $\prod\limits_{
i<\kappa}f(i)/\jbk$ is $\mu$--directed\\
{\em then}\ $\mu\leq\dk$.
\end{enumerate}
\end{proposition}

\Proof (a)\quad Let $\mu=\kappa^+$ and let us assume that $\otp(\kk/\cD)>
\kappa^+$ (the other case is handled similarly). Pick up a function $f\in\kk$
such that 
\[(\forall i<j<\kappa)(i<f(i)<f(j)\ \mbox{ and }f(i)\mbox{ is a regular
cardinal})\]
and $\otp(\prod\limits_{i<\kappa} f(i)/\cD)=\mu$ (remember that $\cD$ is a
normal ultrafilter on $\kappa$). Let $E\subseteq\kappa$ be a
club of $\kappa$ such that  
\[(\forall\delta\in E)(\forall i<\delta)(f(i)<\delta),\]
and let $\{A_\alpha:\alpha<\theta\}\subseteq\cD$ be a family generating $\cD$
and such that $A_\alpha\subseteq E$ (for all $\alpha<\theta$). Choose a
sequence $\bar{g}\subseteq\prod\limits_{i<\kappa}f(i)$ such that $\lh(\bar{g})
=\mu$ and $\langle g_\zeta/\cD:\zeta<\mu\rangle$ is $<_{\cD}$--increasing and
cofinal in $\prod\limits_{i<\kappa}f(i)/\cD$, and $g_\zeta(i)>i$ (for $\zeta<
\mu$ and $i<\kappa$).

For $\alpha<\theta$ and $\zeta<\mu$ we choose a function $h_{\alpha,\zeta}\in
\kk$ such that $h_{\alpha,\zeta}(i)=g_\zeta(\min(A_\alpha\setminus i))$ (for
$i<\kappa$). Let $F=\{h_{\alpha,\zeta}:\alpha<\theta,\ \zeta<\mu\}$, so $|F|
\leq\theta+\mu$. 

Next, for a function $h\in\kk$ define $h^f\in\prod\limits_{i<\kappa}f(i)$ by:
\[h^f(i)=\left\{\begin{array}{ll}
h(i)&\mbox{ if }h(i)<f(i),\\
0   &\mbox{ if }h(i)\geq f(i).
		\end{array}\right.\]
and choose an ordinal $\zeta(h)<\mu$ such that $h^f<_{\cD} g_{\zeta(h)}$. Then
the set $A^h\stackrel{\rm def}{=}\{i<\kappa: h^f(i)<g_{\zeta(h)}(i)\}$ is in
$\cD$. Since the set 
\[A_h\stackrel{\rm def}{=}\{i<\kappa: i\mbox{ is limit and }(\forall j<i)(h(j)
<i)\}\]
is a club of $\kappa$ (so in $\cD$) we may choose $\alpha(h)<\theta$ such that
$A_{\alpha(h)}\subseteq A^h\cap A_h$. 

Suppose now that $G\in [\kk]^{\textstyle\kappa}$, say $G=\{h_\xi:\xi<\kappa
\}$. Take $\zeta<\mu$ such that $\sup\limits_{\xi<\mu}\zeta(h_\xi)<\zeta$ and
let $\alpha<\theta$ be such that 
\[A_\alpha\subseteq\diagun_{\xi<\kappa} A_{\alpha(h_\xi)} \cap \diagun_{\xi
<\kappa}\{i<\kappa: g_{\zeta(h_\xi)}(i)<g_\zeta(i)\}.\]

\begin{claim}
\label{cl1}
If $\xi<i<\kappa$ then $h_\xi(i)\neq h_{\alpha,\zeta}(i)$.
\end{claim}

\noindent {\em Proof of the claim:}\qquad First assume that $\xi<i$, $i\in
A_{\alpha}$. Then, by the choice of $\alpha$, we have $g_{\zeta(h_\xi)}(i)
<g_\zeta(i)$ and $i\in A_{\alpha(h_\xi)}\subseteq A^{h_\xi}$. Consequently,
\[\mbox{either }\ h_\xi(i)\geq f(i)\quad\mbox{or }\ h_\xi(i)=(h_\xi)^f(i)
<g_{\zeta(h_\xi)}(i)<g_\zeta(i)=h_{\alpha,\zeta}(i)<f(i)\]
(and so $h_\xi(i)\neq h_{\alpha,\zeta}(i)$). So suppose now that $\xi<i$,
$i\notin A_\alpha$. Let $j=\min(A_\alpha\setminus i)$. Then $j\in A_{\alpha(
h_\xi)}\subseteq A_{h_\xi}$ and $i<j$, so $h_\xi(i)<j$ and $h_{\alpha,\zeta}
(i)=g_\zeta(j)>j$ . Hence $h_\xi(i)\neq h_{\alpha,\zeta}(i)$.
\medskip

It follows from \ref{cl1} that the family $F$ exemplifies $\geks\leq\theta+
\mu$.
\medskip

\noindent (b)\quad It is similar to (a). Using the assumptions we choose
$f,E,A_\alpha,\bar{g}$ as there and we define $h^f,\zeta(h),\alpha(h)$ (for
$h\in \kk$) in the same manner. Exactly as in \ref{cl1} we show that for each
$h\in\kk$ and $i\in\kk$ we have $h_{\alpha(h),\zeta(h)}(i)\neq h(i)$ (just
consider two cases: $i\in A_{\alpha(h)}$ and $i\notin A_{\alpha(h)}$). 
\medskip

\noindent (c)\quad Similarly. \QED

\stepcounter{section}
\subsection*{\quad 3. Representing functions on the plane}
In this section we answer \cite[Problem 7]{Ci98} showing that it is consistent
that $\con>\aleph_1$ but for every function $f:\bR^2\longrightarrow\bR$ there
exist functions $g_n,h_n:\bR\longrightarrow\bR$, $n<\omega$, such that 
\[f(x,y)=\sum_{n=0}^{\infty} g_n(x)\cdot h_n(y).\]
Let us start with the following technical lemma.

\begin{lemma}
\label{techlem}
Assume ${\bf MA}(\sigma\mbox{-centered})$. Suppose that $B$ is an infinite
subset of $\omega$, $X$ is a set of size $<\con$ and $f,g_n:X\longrightarrow
\bR$ (for $n\in B$) are such that 
\begin{enumerate}
\item[$(\otimes)$] the sets 
\[A_x\stackrel{\rm def}{=}\{n\in B: g_n(x)\neq 0\}\]
for $x\in X$ are infinite almost disjoint.
\end{enumerate}
{\em Then} there is a sequence $\langle b_n:n\in B\rangle$ of rational numbers
such that  
\begin{enumerate}
\item $b_n\neq\frac{1}{(n+1)^2}$ for all $n\in B$, and 
\item $f(x)=\sum\limits_{n\in B} g_n(x)\cdot b_n$\quad for each $x\in X$.
\end{enumerate}
\end{lemma}

\Proof Let $\bQ=\bQ(X,f,B,\langle g_n:n\in B\rangle)$ be the following forcing
notion:

\noindent {\bf a condition in $\bQ$} is a triple $p=(\bar{b}^p,m^p,\sigma^p)=(
\bar{b},m,\sigma)$ such that
\begin{itemize}
\item $m\in\omega$, $\sigma$ is a finite function such that $\dom(\sigma)
\subseteq X$ and $\rng(\sigma)\subseteq B\cap (m+1)$,
\item $\bar{b}=\langle b_n: n\in B\cap m\rangle$ is a sequence of rational
numbers, $b_n\neq\frac{1}{(n+1)^2}$ for $n\in B\cap m$,  
\item for each $x\in\dom(\sigma)$, the sequence
\[\langle |f(x)-\sum\limits_{n\in B\cap k} g_n(x)\cdot b_n|: k\in B\cap
[\sigma(x),m]\rangle\]
is non-increasing;
\end{itemize}
\noindent {\bf the order of $\bQ$} is the natural one:\qquad $p\leq q$\quad if
and only if
\[\bar{b}^p\trianglelefteq\bar{b}^q,\quad m^p\leq m^q\quad\mbox{ and }\quad
\sigma^p\subseteq\sigma^q.\] 

\begin{claim}
\label{cl2}
$\bQ$ is a non-trivial $\sigma$--centered forcing notion.
\end{claim}

\noindent {\em Proof of the claim:}\qquad Consider the space ${}^{\textstyle
X} \omega$ equipped with the product topology of discrete copies of $\omega$.
By Engelking Kar{\l}owicz \cite{EK}, this space is separable (as
$|X|\leq\con$). So let $\{\eta_k:k<\omega\}\subseteq{}^{\textstyle X}\omega$
be a dense subset of ${}^{\textstyle X}\omega$. For $m,k\in\omega$ and a
sequence $\bar{b}=\langle b_n:n\in B\cap m\rangle$ of rationals let
\[Q^{m,\bar{b}}_k\stackrel{\rm def}{=}\{p\in\bQ: m^p=m\ \&\ \bar{b}^p=\bar{b}\
\&\ \sigma^p\subseteq\eta_k\}.\]
Since there are countably many possibilities for $\langle m,\bar{b},k\rangle$
as above and  each member of $\bQ$ belongs to some $Q^{m,\bar{b}}_k$ (remember
the choice of $\eta_k$'s), it is enough to show that the sets $Q^{m,
\bar{b}}_k$ are directed. So let $p_0,\ldots,p_{\ell-1}\in Q^{m,\bar{b}}_k$. 
Then $\bar{b}^{p_i}=\bar{b}$, $m^{p_i}=m$ and $\sigma^{p_i}\subseteq\eta_k$
(for $i<\ell$). Put $q=(\bar{b},m,\bigcup\limits_{i<\ell}\sigma^{p_i})$. It
should be clear that $q\in Q^{m,\bar{b}}_k$ is a condition stronger than all
$p_0,\ldots,p_{\ell-1}$.  
\medskip

Now, for $x\in X$ and a positive rational number $\varepsilon$ let
\[\cI^\varepsilon_x\stackrel{\rm def}{=}\bigg\{p\in\bQ: x\in\dom(\sigma^p)\ \&\
|f(x)-\sum_{n\in B\cap m^p} g_n(x)\cdot b_n|<\varepsilon\bigg\}.\]

\begin{claim}
\label{cl3}
For every $x\in X$ and a rational $\varepsilon>0$ the set $\cI^\varepsilon_x$
is an open dense subset of $\bQ$.
\end{claim}

\noindent {\em Proof of the claim:}\qquad Let $q\in \bQ$ and let $r\in\bQ$ be
defined as follows. If $x\in\dom(\sigma^q)$ then $r=q$, otherwise
\[\bar{b}^r=\bar{b}^q,\quad m^r=m^q\quad\mbox{ and }\quad
\sigma^r=\sigma^q\cup \{(x,m^q)\}.\]
(So $r$ is a condition stronger than $q$ and $x\in\dom(\sigma^r)$.) Use the 
assumption $(\otimes)$ to choose $m^*>m^r$ such that
\[m^*\in B\cap A_x\cap \bigcap\big\{B\setminus A_y:y\in\dom(\sigma^r)\setminus
\{x\}\big\},\]
remember $A_y=\{n\in B: g_n(y)\neq 0\}$. Let 
\[\varepsilon^*=\frac{1}{2}\min\{\varepsilon,|f(x)-\!\sum\limits_{n\in B\cap
m^r}g_n(x)\cdot b_n|\}\]
(so $\varepsilon>\varepsilon^*\geq 0$) and let 
\begin{itemize}
\item $m^p=\min(B\setminus (m^*+1))$, $\sigma^p=\sigma^r$,
\item $b^p_n=b^r_n$ if $n\in B\cap m^r$, $b^p_n=0$ if $n\in B\cap [m^r,m^*)$
and $b^p_{m^*}\neq\frac{1}{(m^*+1)^2}$ be a rational number such that
\[f(x)-\!\sum_{n\in B\cap m^*} g_n(x)\cdot b^p_n-\varepsilon^*\leq g_{m^*}(x)
\cdot b^p_{m^*}\leq f(x)-\!\sum_{n\in B\cap m^*} g_n(x)\cdot b^p_n+
\varepsilon^*\]
(clearly the choice is possible as $g_{m^*}(x)\neq 0$; if $\varepsilon^*=0$
then $b^p_{m^*}=0$). 
\end{itemize}
One easily checks now that the above choice defines a condition
$p\in\cI^\varepsilon_x$ stronger than $r$.
\medskip

It follows from \ref{cl2}, \ref{cl3} that we may use ${\bf
MA}(\sigma\mbox{-centered})$ to find a directed set $G\subseteq \bQ$ such that
$G\cap \cI^\varepsilon_x\neq\emptyset$ for each $x\in X$ and a positive
rational $\varepsilon$. Let $\bar{b}=\bigcup\{\bar{b}^p: p\in G\}$. It should
be clear that the sequence $\bar{b}$ is as required. \QED

\begin{definition}
\label{approx}
Let $f:\bR^2\longrightarrow\bR$ . An {\em $f$--approximation} is a tuple
$p=(X^p,\bar{g}^p_0,\bar{g}^p_1,\cD^p)=(X,\bar{g}_0,\bar{g}_1,\cD)$ such
that  
\begin{enumerate}
\item[(a)] $X\subseteq\bR$ ,
\item[(b)] $\bar{g}_\ell=\langle g_{\ell,n}:n<\omega\rangle$,\quad $g_{\ell,
n}: X\longrightarrow\bR$\quad (for $\ell<2$, $n<\omega$),

\noindent for $\ell<2$, $x\in X$ let $\bar{a}_{\ell,x}=\langle g_{\ell,n}(x):
n<\omega\rangle$,
\item[(c)] $(\forall x,y\in X)(f(x,y)=\sum\limits_{n=0}^\infty g_{0,n}(x)
\cdot g_{1,n}(x))$, 
\item[(d)] $\cD$ is a filter on $\omega$ including all co-finite subsets of
$\omega$ and generated by $\leq |X|+\aleph_0$ sets,
\item[(e)] if $x\in X$, $\ell<2$ then
\[\begin{array}{ll}
\ &\big\{n<\omega: g_{\ell,n}(x)\in\{0,\frac{1}{(n+1)^2}\}\big\}\in\cD,\qquad
\mbox{and}\\
A_{\bar{a}_{\ell,x}}=A^p_{\bar{a}_{\ell,x}}\stackrel{\rm def}{=}&\big\{n<
\omega:g_{\ell,n}(x)=\frac{1}{(n+1)^2}\big\}\neq\emptyset\mod\cD,
  \end{array}\]
\item[(f)] no finite union of sets $A_{\bar{a}_{\ell,x}}$ (for $\ell<2$, $x\in
X$) is in $\cD$,
\item[(g)] if $(\ell_1,x_1)\neq (\ell_2,x_2)$, $\ell_1,\ell_2<2$, $x_1,x_2\in
X$ then $\bar{a}_{\ell_1,x_1}\neq \bar{a}_{\ell_2,x_2}$ and the intersection
$A^p_{\bar{a}_{\ell_1,x_1}}\cap A^p_{\bar{a}_{\ell_2,x_2}}$ is finite.
\end{enumerate}
The set $\APf$ of all approximations carries a natural partial order:\\
for $f$--approximations $p,q$ we let\qquad $p\leq q$\quad if and only if
\[X^p\subseteq X^q,\quad \cD^p\subseteq\cD^q\quad\mbox{ and }\quad g^p_{\ell,
n}\subseteq g^q_{\ell,n}\quad\mbox{ (for $\ell<2$ and $n<\omega$).}\] 
\end{definition}

\begin{theorem}
\label{MAthm}
Assume ${\bf MA}(\sigma\mbox{-centered})$. Let $f:\bR^2\longrightarrow\bR$ . 
Suppose that $p\in\APf$ is such that $|X^p|<\con$ and $r^*\in\bR\setminus
X^p$. {\em Then\/} there is $q\in\APf$ such that  
\[p\leq q\quad\mbox{ and }\quad X^q=X^p\cup\{r^*\}.\]
\end{theorem}

\Proof First choose pairwise disjoint infinite subsets $B_0,B_1,B_2$ of
$\omega\setminus\{0\}$ such that for $m<3$:
\begin{enumerate}
\item[$(\alpha)$] $(\forall B\in\cD^p)(B_m\subseteq^* B)$,
\item[$(\beta)$]  $(\forall\ell<2)(\forall x\in X^p)(|B_m\cap A_{\bar{a}^p_{
\ell,x}}|=\aleph_0)$,
\item[$(\gamma)$] no finite union of sets $A_{\bar{a}_{\ell,x}}$ (for
$\ell<2$, $x\in X$) almost includes $B_m$.
\end{enumerate}
(There are such sets by ${\bf MA}(\sigma\mbox{-centered})$; remember
\ref{approx}(e),(f).) Next choose disjoint infinite subsets $B^0_0,B^1_0,
B^2_0$ of $B_0$ such that for $k<3$
\begin{enumerate}
\item[$(\delta)$] $(\forall\ell<2)(\forall x\in X^p)(|B^k_0\cap A_{\bar{a}^p_{
\ell,x}}|<\aleph_0)$.
\end{enumerate}
(Again, easily possible by our assumptions and \ref{approx}(g) and $(\gamma)$
above.) 

Now we start defining an $f$--approximation $q$. We let
\begin{itemize}
\item $X^q=X^p\cup\{r^*\}$,
\item $\cD^q$ be the filter generated by $\cD^p\cup\{B_0\}$,
\item $g^q_{\ell,n}(x)=g^p_{\ell, n}(x)$\quad for $x\in X^p$, $\ell<2$ and
$n<\omega$, 
\item if $n\in\omega\setminus B_1$ then
\[g^q_{0,n}(r^*)=\left\{\begin{array}{ll}
1&\mbox{ if }n=0,\\
\frac{1}{(n+1)^2}&\mbox{ if }n\in B^0_0,\\
0&\mbox{ if }n\in\omega\setminus (B_1\cup B^0_0\cup\{0\});
		       \end{array}\right.\]
and if $n\in\omega\setminus B_2$ then
\[g^q_{1,n}(r^*)=\left\{\begin{array}{ll}
f(r^*,r^*)&\mbox{ if }n=0,\\
\frac{1}{(n+1)^2}&\mbox{ if }n\in B^1_0,\\
0&\mbox{ if }n\in\omega\setminus (B_2\cup B^1_0\cup\{0\}).
		       \end{array}\right.\]
\end{itemize}
Now we want to define $g^q_{0,n}(r^*)$, $g^q_{1,n}(r^*)$ for other $n$, but we
have to be careful with that to ensure that the clause \ref{approx}(c) is
satisfied. It should be clear at the moment that we do not have to worry
anymore about that clause if $x,y\in X^p$ or $x=y=r^*$ (for the last case
inspect the definition above and the choice of $B_m,B^0_0,B^1_0$). So now we
use \ref{techlem} to  finish the definition. First note that for each $x\in X$
the sets
\[\{n\in\omega\setminus B_1\!: g^q_{0,n}(r^*)\cdot g^p_{1,n}(x)\neq 0\}\ \mbox{
and }\ \{n\in\omega\setminus B_2\!: g^p_{0,n}(x)\cdot g^q_{1,n}(r^*)\neq 0\}\]
are finite (remember clauses $(\alpha)$ and $(\delta)$). Apply \ref{techlem}
to the set $B_1$, functions $g^p_{1,n}$ (for $n\in B_1$) and the mapping  
\[x\mapsto f(r^*,x)-\sum_{n\in\omega\setminus B_1} g^q_{0,n}(r^*)\cdot g^p_{1,
n}(x)\]
(note that the sum is actually finite) to find $g^q_{0,n}(r^*)$ (for $n\in
B_1$) such that $g^q_{0,n}(r^*)\neq\frac{1}{(n+1)^2}$ and for each $x\in X$
\[f(r^*,x)-\sum_{n\in\omega\setminus B_1} g^q_{0,n}(r^*)\cdot g^p_{1,n}(x)=
\sum_{n\in B_1} g^q_{0,n}(r^*)\cdot g^p_{1,n}(x).\]
Next use \ref{techlem} for $B_2$, $g^p_{0,n}$ (for $n\in B_2$) and the mapping  
\[x\mapsto f(x,r^*)-\sum_{n\in\omega\setminus B_2} g^p_{0,n}(x)\cdot g^q_{1,
n}(r^*)\]
to choose $g^q_{1,n}(r^*)$ (for $n\in B_2$) such that $g^q_{1,n}(r^*)\neq
\frac{1}{(n+1)^2}$ and for $x\in X$
\[f(x,r^*)-\sum_{n\in\omega\setminus B_2} g^p_{0,n}(x)\cdot g^q_{1,n}(r^*)=
\sum_{n\in B_2} g^p_{0,n}(x)\cdot g^q_{1,n}(r^*).\]
This finishes the definition of $g^q_{\ell,n}(x)$ for $\ell<2$, $x\in X^q$ and
$n<\omega$. Checking that $(X^q,\bar{g}^q_0,\bar{g}^q_1,\cD^q)\in\APf$ is as
required is straightforward. \QED
\medskip

Since $\leq$--increasing sequences of $f$--approximations have (natural) upper
bounds we may use \ref{MAthm} to prove inductively the following.

\begin{conclusion}
\label{condwa}
Assume ${\bf MA}(\sigma\mbox{-centered})$. Then for every function $f:\bR^2
\longrightarrow\bR$ there are functions $g_n,h_n:\bR\longrightarrow\bR$,
$n<\omega$, such that  
\[f(x,y)=\sum_{n=0}^{\infty} g_n(x)\cdot h_n(y).\]
\end{conclusion}

\begin{remark}
Regarding the assumptions of \ref{condwa}, remember that by Bell \cite{B}
${\bf MA}(\sigma\mbox{-centered})$ is equivalent to ${\frak p}=\con$.
\end{remark}

Let us finish this section with the following ``negative'' result.

\begin{proposition}
Let $\bP$ be the forcing notion for adding $\aleph_2$ Cohen reals. Then, in
$\bV^{\bP}$, there is a function $f:\bR^2\longrightarrow\bR$ such that there
are {\bf no} functions $g_n,h_n:\bR\longrightarrow\bR$ satisfying
\[(\forall x,y\in\bR)(f(x,y)=\sum_{n=0}^\infty g_n(x)\cdot h_n(y)).\]
\end{proposition}

\Proof Since we may break $\bP$ into two steps each adding $\aleph_2$ Cohen
reals, we may assume that $\bV\models\neg {\bf CH}$. So fix a sequence
$\langle\eta_i:i<\aleph_2\rangle$ of pairwise distinct real numbers. Then
$\bP$ may be interpreted as the partial order of all finite functions $p$ such
that $\dom(p)\subseteq\{(\eta_i,\eta_j): i,j<\aleph_2\}$ and $\rng(p)\subseteq
2$ ordered by the inclusion. For a set $A\subseteq\{(\eta_i,\eta_j): i,j<
\aleph_2\}$ let $\bP_A=\{p\in\bP:\dom(p)\subseteq A\}$ (so
$\bP_A\lesdot\bP$). 

Let $\name{f}$ be a $\bP$--name for a function from $\bR^2$ to $\bR$ such that
$\forces_{\bP}\bigcup\{p:p\in\name{G}_{\bP}\}\subseteq\name{f}$. Suppose that
$\name{g}_n,\name{h}_n$ (for $n<\omega$) are $\bP$--names for functions from
$\bR$ to $\bR$.  

\begin{claim}
\[\forces_{\bP}\mbox{`` }(\exists i<j<\aleph_2)(\name{f}(\eta_i,\eta_j)\neq
\sum_{n=0}^\infty \name{g}_n(\eta_i)\cdot\name{h}_n(\eta_j))\mbox{ ''.}\]
\end{claim}

\noindent {\em Proof of the claim:}\qquad Let $q\in \bP$. For each
$i<\aleph_2$ fix a countable subset $A_i$ of $\{(\eta_\xi,\eta_\zeta): \xi,
\zeta<\aleph_2\}$ such that $\dom(q)\subseteq A_i$ and for some
$\bP_{A_i}$--names $\name{r}_n,\name{s}_n$ (for $n<\omega$) we have
$\forces_{\bP}$`` $\name{g}_n(\eta_i)=\name{r}_n$ and $\name{h}_n(\eta_i)=
\name{s}_n$ ''. Let 
\[B_i\stackrel{\rm def}{=}\{\xi: (\exists\zeta<\aleph_2)((\eta_\xi,
\eta_\zeta)\in A_i\ \mbox{ or }\ (\eta_\zeta,\eta_\xi)\in A_i)\}\]
(clearly each $B_i$ is countable). Plainly, for $i\in S^2_1\stackrel{\rm
def}{=}\{\delta<\aleph_2:\cf(\delta)=\aleph_1\}$ we have $\sup(B_i\cap i)<i$
and hence for some $j<\aleph_2$ the set $S=\{i\in S^2_1: \sup(B_i\cap i)=j\}$
is stationary. Choose $i_0<i_1$ from $S$ such that $\sup(B_{i_0})<i_1$. Let
\[Y=\{(\eta_\xi,\eta_\zeta): \{\xi,\zeta\}\subseteq B_{i_0}\ \mbox{ or }\
\{\xi,\zeta\}\subseteq B_{i_1}\}.\]
Note that $(\eta_{i_0},\eta_{i_1})\notin Y$. Since $\name{g}_n(\eta_{i_0}),
\name{h}_n(\eta_{i_1})$ are (essentially)
$\bP_Y$--names and $q\in\bP_Y$, we find a condition $p\in\bP_Y$ stronger than
$q$ and deciding the statement ``$\sum\limits_{n=0}^\infty \name{g}_n(
\eta_{i_0})\cdot\name{h}_n(\eta_{i_1})\leq \frac{1}{2}$''. Let $r\in\bP$ be a
condition stronger than $p$ such that $(\eta_{i_0},\eta_{i_1})\in\dom(r)$ and 
\[r(\eta_{i_0},\eta_{i_1})=\left\{\begin{array}{ll}
1&\mbox{ if }p\forces_{\bP_Y}\mbox{`` }\sum\limits_{n=0}^\infty\name{g}_n(
\eta_{i_0})\cdot\name{h}_n(\eta_{i_1})\leq \frac{1}{2}\mbox{ '',}\\
0&\mbox{ otherwise.}
				  \end{array}\right.\]
Then 
\[r\forces_{\bP}\name{f}(\eta_{i_0},\eta_{i_1})\neq\sum\limits_{n=0}^\infty
\name{g}_n(\eta_{i_0})\cdot\name{h}_n(\eta_{i_1}),\]
finishing the proof. \QED

\stepcounter{section}
\subsection*{\quad 4. Countably continuous functions}
Our aim here is to show that, consistently, {\bf CH} fails but every $f:\bR
\longrightarrow\bR$ satisfying 
\[(\forall U\in [\bR]^{\textstyle \aleph_1})(\exists U^*\in [U]^{\textstyle
\aleph_1})(f\rest U^*\mbox{ is continuous\/})\]
is countably continuous (see Definition \ref{countcon} below). This
answers negatively \cite[Problem 4]{Ci98}.  

\begin{definition}
\label{countcon}
A function $f:\bR\longrightarrow\bR$ is {\em countably continuous} if there is
a partition $\langle X_n:n<\omega\rangle$ of $\bR$ such that the restriction
of $f$ to any $X_n$ is continuous. 
\end{definition}

\begin{theorem}
\label{pr4thm}
It is consistent with $\neg {\bf CH}$ that every function $f:\bR
\longrightarrow\bR$ such that
\begin{enumerate}
\item[$(\oplus)_f$]\quad $(\forall U\in [\bR]^{\textstyle \aleph_1})(\exists
U^*\in [U]^{\textstyle\aleph_1})(f\rest U^*\mbox{ is continuous\/})$
\end{enumerate}
is countably continuous.
\end{theorem}

\Proof Start with $\bV\models {\bf CH}$ and let $\lambda>\aleph_1$ be a
cardinal such that $\lambda^{\aleph_0}=\lambda$. Let $\bP_\lambda$ be a
forcing notion for adding $\lambda$ many Cohen reals. So $\bP_\lambda$ can be
represented as the set of all finite partial functions $p:\dom(p)
\longrightarrow 2$, $\dom(p)\subseteq\lambda$, ordered by the inclusion. 

For a set $A\subseteq \lambda$ let $\bP_A=\{p\in\bP_\lambda:\dom(p)\subseteq
A\}$. Then $\bP_A\lesdot\bP_\lambda$. Plainly, $\bP_\lambda$ is a ccc forcing
notion and $\forces_{\bP_\lambda}\con=\lambda>\aleph_1$. 

We are going to show that in $\bV^{\bP_\lambda}$, every real function $f$
satisfying $(\oplus)_f$ is countably continuous. To this end suppose that
$\name{f}$ is a $\bP_\lambda$--name for a function from $\bR$ into $\bR$ such
that  
\[\forces_{\bP_\lambda}\mbox{`` $\name{f}$ is not countably continuous ''}.\]
By induction on $\alpha<\omega_1$ choose an increasing continuous sequence
$\langle A_\alpha:\alpha<\omega_1\rangle$ such that for each
$\alpha<\omega_1$: 
\begin{enumerate}
\item $A_\alpha\in [\lambda]^{\textstyle \aleph_1}$; 
\item if $\name{\eta}$ is a $\bP_{A_\alpha}$--name for a real then
$\name{f}(\name{\eta})$ is a $\bP_{A_{\alpha+1}}$--name;
\item if $\name{\bar{h}}=\langle\name{h}_n: n<\omega\rangle$ is a
$\bP_{A_\alpha}$--name for an $\omega$--sequence of partial real functions
such that $\dom(\name{h}_n)$ is a Borel set and $\name{h}_n$ is continuous on
its domain (for $n<\omega$) then each $\name{h}_n$ is a
$\bP_{A_{\alpha+1}}$--name and there is a $\bP_{A_{\alpha+1}}$--name
$\name{\eta}$ for a real such that 
\[\forces_{\bP_{A_{\alpha+1}}}\mbox{`` }(\forall n<\omega)(\name{f}(
\name{\eta})\neq\name{h}_n(\name{\eta}))\mbox{ ''}.\]
\end{enumerate}
(There are no problems with carrying out the construction.) Let $A=\bigcup 
\limits_{\alpha<\omega_1} A_\alpha$, so $A\in [\lambda]^{\textstyle\aleph_1}$
and $\forces_{\bP_A}{\bf CH}$. It should be clear that, by (2) above, we have
a $\bP_A$--name $\name{f}^A$ such that $\forces_{\bP_\lambda}$`` $\name{f}^A=
\name{f}\rest \bR\cap\bV^{\bP_A}$ ''. Moreover, by (3), we know that
\[\forces_{\bP_A}\mbox{`` $\name{f}^A$ is not countably continuous ''.}\]
Now, using \cite[3.11]{Ci98}, we conclude that (remember that we have {\bf CH}
in $\bV^{\bP_A}$)
\[\forces_{\bP_A}\mbox{`` }(\exists U\in [\bR]^{\textstyle \aleph_1})(\forall
U^*\in [U]^{\textstyle\aleph_1})(\name{f}^A\rest U^*\mbox{ is not continuous
\/})\mbox{ ''.}\]
Let $G\subseteq\bP_A$ be a generic filter over $\bV$. Work in $\bV[G]$. Let
$U\in [\bR]^{\textstyle\aleph_1}$ be such that $\name{f}^A[G]\rest U^*$ is not
continuous for any uncountable $U^*\subseteq U$. We want to show that this
property of the function $\name{f}^A[G]$ and the set $U$ is preserved by the
quotient forcing $\bP_\lambda/\bP_A$ (which is isomorphic to $\bP_{\lambda
\setminus A}$, of course). So suppose that $p\in\bP_\lambda/\bP_A$ is such
that 
\[p\forces_{\bP_\lambda/\bP_A}\mbox{`` }(\exists U^*\in [U]^{\textstyle
\aleph_1})(\name{f}^A[G]\rest U^*\mbox{ is continuous\/})\mbox{ ''.}\]  
Every continuous function on a set $U^*\subseteq\bR$ can be extended to a
continuous function on a $\Pi^0_2$--set. Now, both $\Pi^0_2$--sets and
continuous functions on them are coded by reals. Consequently we find a
countable set $B\subseteq\lambda\setminus A$ such that $\dom(p)\subseteq B$
and for some $\bP_B$--names $\name{W}, \name{h}$ we have
\[\begin{array}{ll}
p\forces_{\bP_\lambda/\bP_A}\mbox{``}&\name{W}\mbox{ is a $\Pi^0_2$--subset
of $\bR$ , }\ \name{h}:\name{W}\longrightarrow\bR\mbox{ is continuous and }\\
\ &(\exists^{\aleph_1}\eta\in U)(\eta\in\name{W}\ \ \&\ \ \name{f}^A
[G](\eta)=\name{h}(\eta))\mbox{ ''.}
  \end{array}\]
The property stated above is absolute from $\bV^{\bP_\lambda}$ to $\bV^{\bP_{A
\cup B}}$, so the condition $p$ forces the respective sentence in $\bP_B$. Now, 
the forcing notion $\bP_B$ is countable so it has the property that every
uncountable set of ordinals in the extension contains an uncountable subset
from the ground model. Consequently, we find (still in $\bV[G]$) an uncountable
set $U_0\subseteq U$ such that
\[\begin{array}{ll}
p\forces_{\bP_B}\mbox{``}&\name{W}\mbox{ is a $\Pi^0_2$--subset
of $\bR$, }\ \name{h}:\name{W}\longrightarrow\bR\mbox{ is continuous and }\\
\ &(\forall \eta\in U_0)(\eta\in\name{W}\ \ \&\ \ \name{f}^A[G](\eta)=\name{h}
(\eta))\mbox{ ''.}
  \end{array}\]
Thus $p\forces_{\bP_B}$ ``$\name{f}^A[G]\rest U_0$ is continuous '', and
hence easily this statement has to hold in $\bV[G]$ already, a contradiction.
\QED

\end{document}